\documentclass[12pt]{article}

\usepackage{amssymb}
\usepackage{amsmath}
\usepackage{authblk}
\usepackage{graphicx}
\usepackage{subcaption}
\usepackage{listings}
\DeclareMathOperator{\arccosh}{arccosh}
\newcommand{\RR}{\mathbb R}

\begin{document}

\title{Second order stabilized two-step Runge-Kutta methods}

\author[1]{Andrew Moisa}
\author[2]{Boris Faleichik}
\affil[1]{Kreo Software Ltd., 71-75 Shelton Street, Covent Garden, London, WC2H~9JQ, United~Kingdom}
\affil[2]{Department of Computational Mathematics, Belarusian State University, 4~Nezavisimosti~Avenue, Minsk, 220030, Belarus}

\maketitle

\begin{abstract}
Stabilized methods (also called Chebyshev methods) are explicit methods with extended stability domains along the negative real axis. These methods are intended for large mildly stiff problems, originating mainly from parabolic PDEs. In this paper we present explicit two-step Runge-Kutta methods, which have an increased stability interval in comparison with one-step methods (up to 2.5 times). Also, we perform some numerical experiments to confirm the accuracy and stability of this methods.
\end{abstract}



\section{Introduction}
\label{}

Up to now, there exist the following types of stabilized explicit methods:
\begin{itemize}
	\item One-step Runge-Kutta methods, which use internal stages to increase the stability interval (see \cite[pp. 31-36]{hairer2}, \cite{vdhouwen}, \cite{abdulle2}, \cite{lebedev}). These methods may require a sufficiently large number of function evaluations to achieve the required length of the stability interval.
	\item Multistep Adams-types methods, which use previously calculated function evaluations to increase the stability interval (see \cite{repnikov}). These methods require only one function evaluation per step, but have other significant disadvantages inherent in multi-step methods.
\end{itemize}

In this paper we present two-step methods, which use internal stages to increase the stability interval. These methods require about 1.5 times less stages to achieve the same stability as one-step.

The paper is organized as follows. In Sections \ref{stabintopt} and \ref{damping} we describe the way to get an optimal stability region for the considered methods. In Section \ref{errstab} we calculate error constants and length of stability interval for them. In Section \ref{ortpolys} we explain how to construct two-step methods that use the recurrence relation. Section \ref{numexps} contains numerical experiments to test the properties of the methods obtained.

\section{Stability interval optimization}
\label{stabintopt}

A two-step $s$-stage explicit Runge-Kutta method for the numerical integration of the ODE system
\begin{equation}\label{eq:ode}
	y'=f(t,y), \quad y(x_0) = y_0\in\RR^n, \quad y: \RR \to \RR^n, \quad f:\RR\times \RR^n \to \RR^n
\end{equation}
with preconsistency condition \cite[(9.23)]{hairer2} has the form
\begin{equation}\label{eq:inmethod}
	\begin{aligned}
		v_1 &= \tilde{a}_1 y_n + (1-\tilde{a}_1) y_{n-1}, \\
		v_2 &= \tilde{a}_2 y_n + (1-\tilde{a}_2) y_{n-1} + h 	\tilde{b}_{21} f(x_n + c_1 h, v_1), \\
		\vdots \\
		v_s &= \tilde{a}_s y_n + (1-\tilde{a}_s) y_{n-1} \: + \\
		& \qquad \quad + h \left( \tilde{b}_{s,1} f(x_n + c_1 h, v_1) + \dots + \tilde{b}_{s,s-1} f(x_n + c_{s-1} h, v_{s-1}) \right), \\
		y_{n+1} &= a y_n + (1-a) y_{n-1} \: + \\
		& \qquad \quad + h \left( b_1 f(x_n + c_1 h, v_1) + \dots + b_s f(x_n + c_s h, v_s) \right)
	\end{aligned}
\end{equation}
\cite[p. 362]{hairer2}. Coefficients $c_j$ can be found by the formula
\begin{equation}\label{eq:inc}
	c_j = \tilde{a}_j - 1 + \sum_{k=1}^{j-1} \tilde{b}_{jk}
\end{equation}
\cite[p. 443]{hairer1}.

Characteristic equation for \eqref{eq:inmethod} has the form
\begin{equation}\label{eq:chareq}
	\begin{aligned}
		&\zeta^2 - R_s^1(\mu) \zeta - R_s^0(\mu) = 0, \\
		R_s^1(\mu) = \sum_{j=0}^{s} r_j^1 \mu^j &= a + \mu \sum_{j} b_j \tilde{a}_j + \mu^2 \sum_{j,k} b_j \tilde{b}_{jk} \tilde{a}_k + \\
		&\qquad\qquad\qquad\qquad + \cdots + \mu^s b_s \left( \prod_{j} \tilde{b}_{j,j-1} \right) \tilde{a}_1, \\
		R_s^0(\mu) = \sum_{j=0}^{s} r_j^0 \mu^j &= (1-a) + \mu \sum_{j} b_j (1-\tilde{a}_j) + \mu^2 \sum_{j,k} b_j \tilde{b}_{jk} (1-\tilde{a}_k) + \\ 
		&\qquad\qquad\qquad\qquad + \cdots + \mu^s b_s \left( \prod_{j} \tilde{b}_{j,j-1} \right) (1-\tilde{a}_1)
	\end{aligned}
\end{equation}
(for comparison with one-step Runge-Kutta methods see \cite[p. 16]{hairer2}).

Order conditions can be written as
\begin{equation}\label{eq:ordercond}
	\begin{cases}
		r_1^0 + r_1^1 + a = 2, \\
		r_2^0 + r_2^1 + r_1^1 + \displaystyle\frac{a}{2} = 2.
	\end{cases}
\end{equation}

Our task is to find polynomials $R^0$ and $R^1$ such that the corresponding stability interval is as large as possible and their coefficients satisfy the posed order conditions \eqref{eq:ordercond}. 

Let's first look at the quadratic equation \eqref{eq:chareq}. For its roots $-1 \le |\zeta_1|, |\zeta_2| \le 1$ polynomial $R^0$ must be inside the segment $[-1, 1]$ and $R^1$ must be inside $[-1 + R^0, 1 - R^0]$. As in the case of one-step methods, the best option for $R^0$ is the Chebyshev polynomial:
\begin{equation}\label{eq:initialR0}
	R_s^0(\mu) = \pm T_s\left( 1 + \frac{\mu}{l} \right).
\end{equation}
This polynomial remains between $-1$ and $+1$ on the largest possible interval $[-2l, 0]$.

\begin{figure}
	\begin{subfigure}{0.48\textwidth}
		\includegraphics[width=\textwidth]{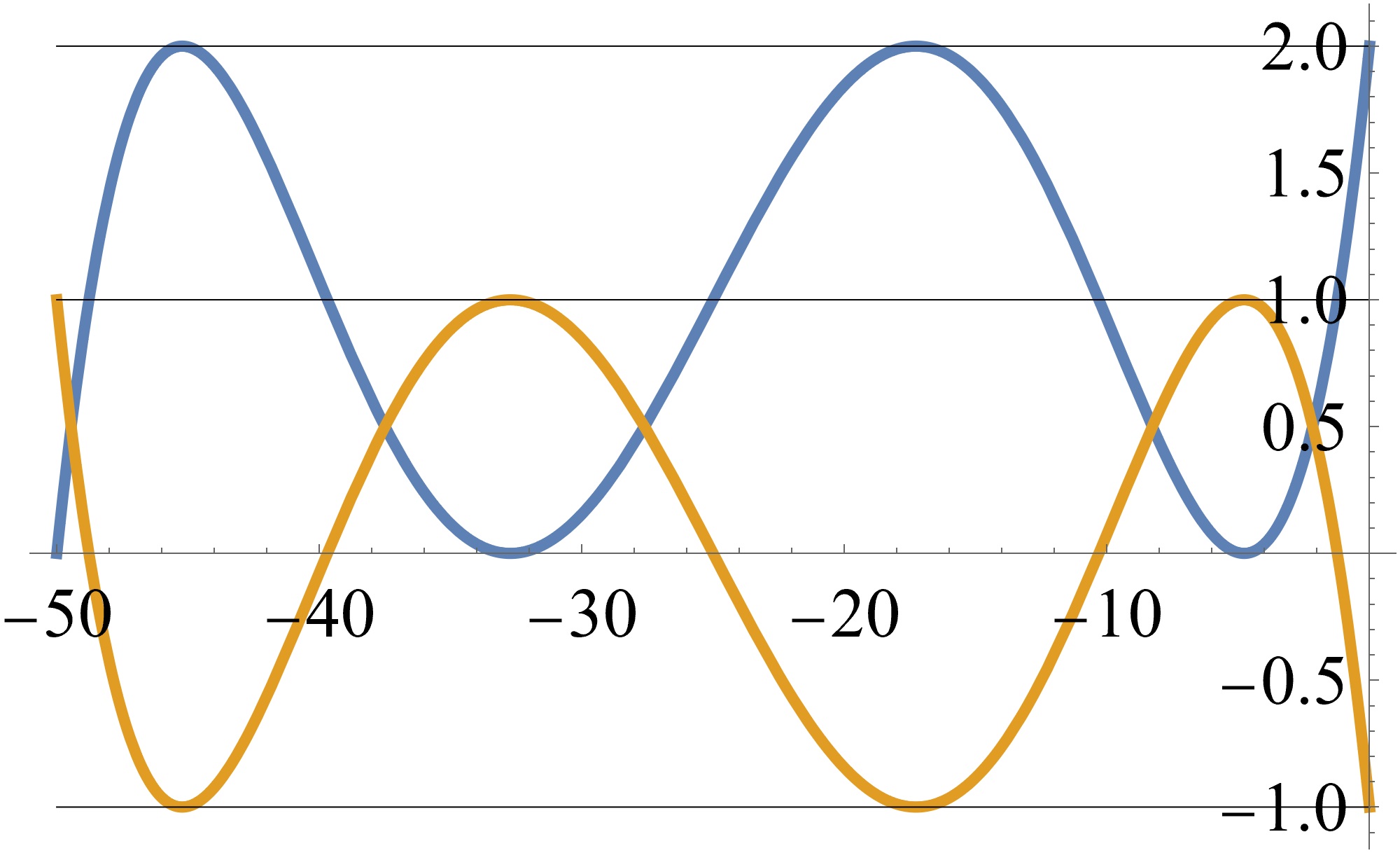}
	\end{subfigure}
	\hfill
	\begin{subfigure}{0.48\textwidth}
		\includegraphics[width=\textwidth]{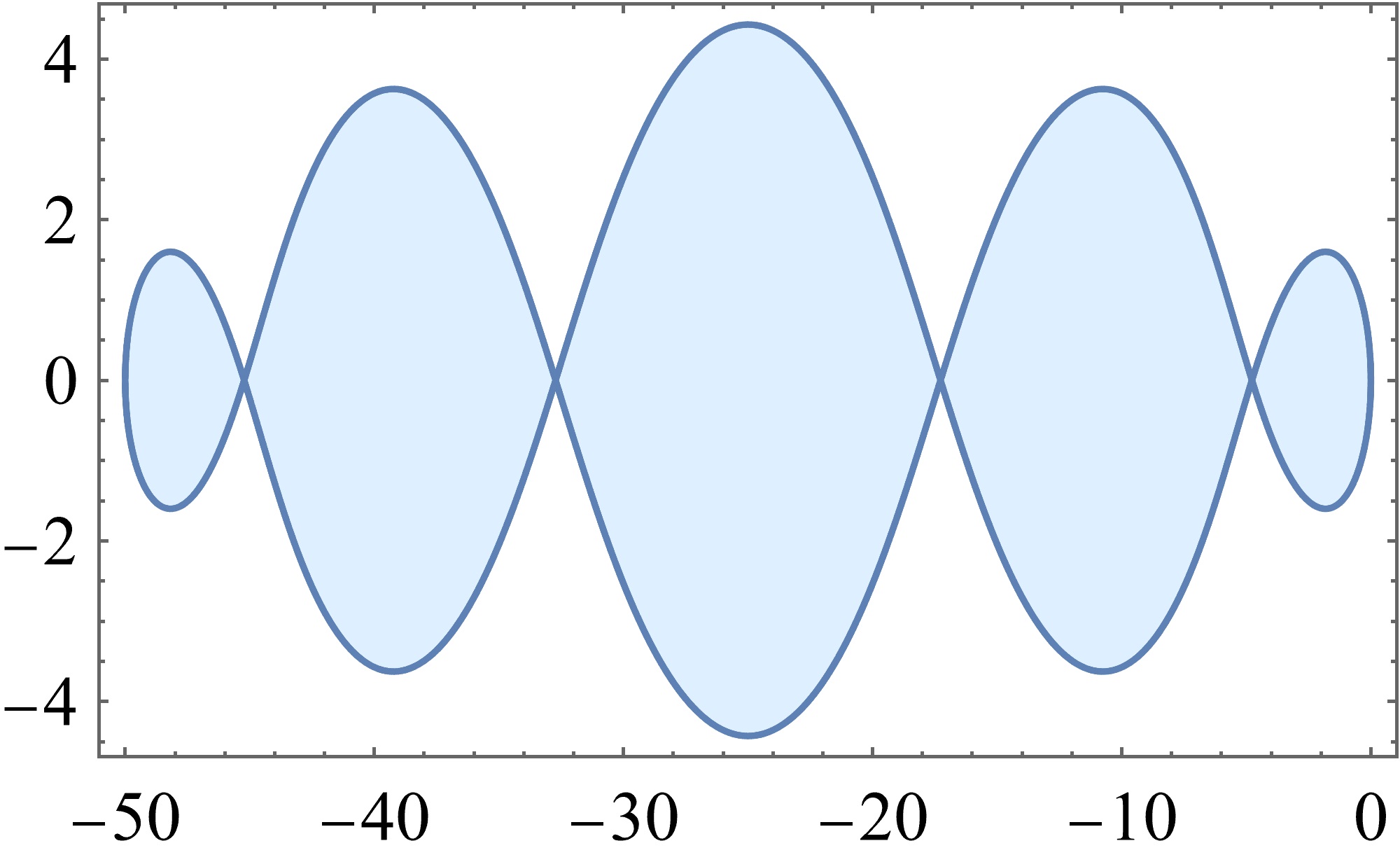}
	\end{subfigure}
	\caption{Shifted Chebyshev polynomials and their stability domain $(s = 5)$}
	\label{fig1}
\end{figure}

Notice, that quadratic equation $x^2 - (1 + c) x + c = 0$ has the roots $x_1 = 1$ and $x_2 = c$. Thus, all our conditions are satisfied by the polynomials
\begin{equation}\label{eq:initialR1R0}
	R_s^1(\mu) = 1 + T_s\left( 1 + \frac{\mu}{s^2} \right), \quad
	R_s^0(\mu) = -T_s\left( 1 + \frac{\mu}{s^2} \right)
\end{equation}
(see Figure~\ref{fig1}).

Note, that the points where $R^0 = \pm1$ are not included in the stability domain. Now it is not important, because damping procedure is needed anyway.

\section{Damping}
\label{damping}

The damping of the pair of polynomials \eqref{eq:initialR1R0} will be carried out in two stages: first, we will produce an "ideal" damping with a possible violation of the order conditions \eqref{eq:ordercond}, and after that we will achieve the observance of the order conditions for the resulting pair.
As always choose $0 < \varepsilon \ll 1, \eta = 1 - \varepsilon$. Polynomials (7) should be changed in such a way that the roots of equation (4) belong to the segment $[-\eta, \eta]$. It is easy to see that
\begin{equation}\label{eq:indaR1R0}
	R_s^1(\mu) = \eta \left(1 + T_s\left( 1 + \frac{\mu}{s^2} \right) \right), \quad
	R_s^0(\mu) = -\eta^2 T_s\left( 1 + \frac{\mu}{s^2} \right)
\end{equation}
satisfy this condition.

It remains to satisfy the order conditions \eqref{eq:ordercond}. Let’s consider the pair
\begin{equation}\label{eq:dampedR1R0}
	R_s^1(\mu) = \alpha \left(1 + T_s\left( \omega + \beta \frac{\mu}{s^2} \right) \right), \quad
	R_s^0(\mu) = -\eta^2 T_s\left( \omega + \beta \frac{\mu}{s^2} \right).
\end{equation}

\begin{figure}
	\begin{subfigure}{0.48\textwidth}
		\includegraphics[width=\textwidth]{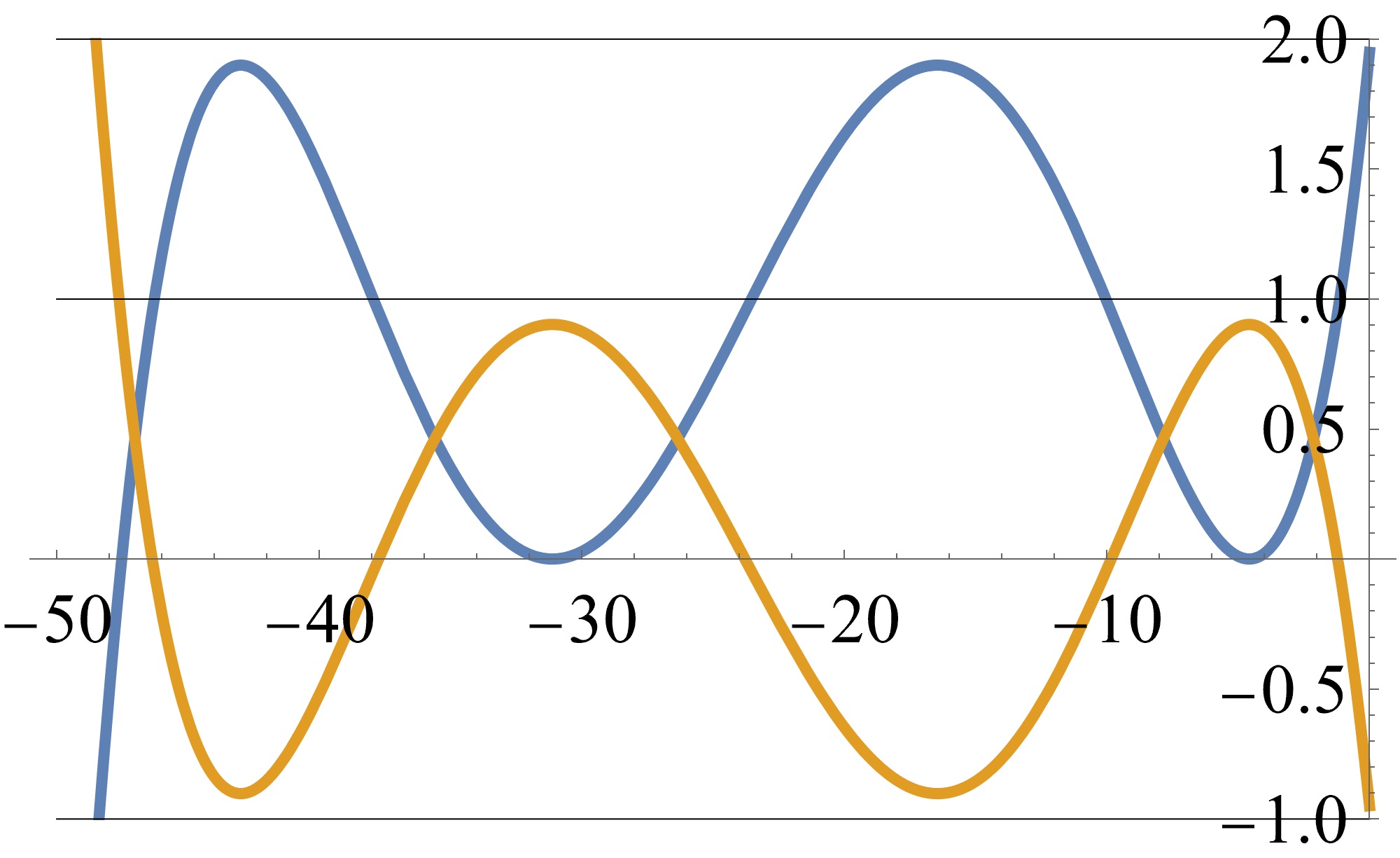}
	\end{subfigure}
	\hfill
	\begin{subfigure}{0.48\textwidth}
		\includegraphics[width=\textwidth]{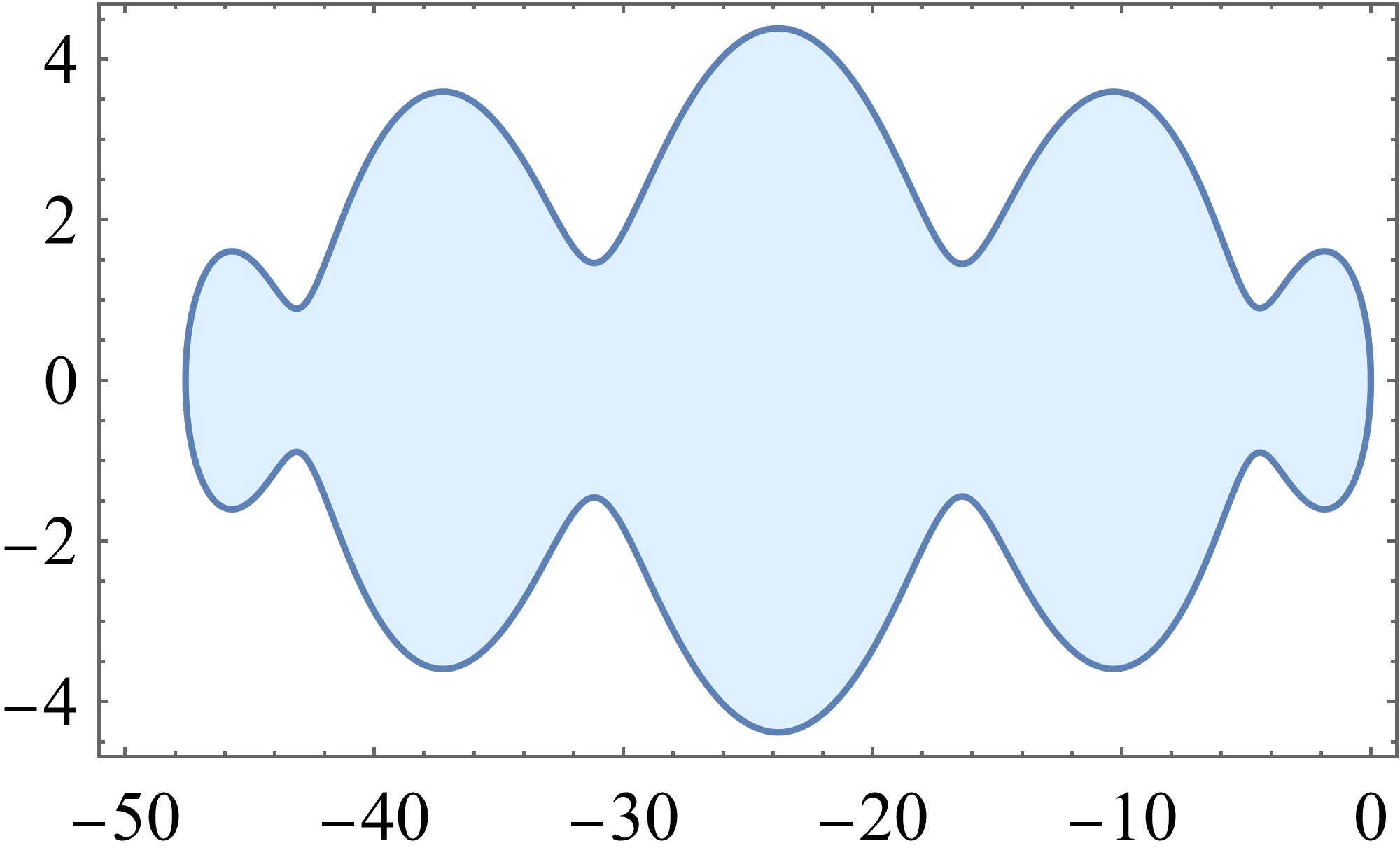}
	\end{subfigure}
	\caption{Shifted Chebyshev polynomials (damping) and their stability domain $(s = 5, \varepsilon = 0.05)$}
	\label{fig2}
\end{figure}

If we write the order conditions \eqref{eq:ordercond} and condition on free members \linebreak $\left( r_0^1 + r_0^0 = 1 \right)$ for it, we obtain a system of 3 nonlinear equations of variables $\{\alpha, \omega, \beta\}$:
\begin{equation}\label{eq:nonlinear}
	\begin{cases}
		\alpha + (\alpha - \eta^2) T_s(\omega) - 1 = 0, \\
		\displaystyle
		\alpha + \alpha T_s(\omega) + \frac{\beta(\alpha-\eta^2)}{s^2}T_s'(\omega) - 2 = 0, \\
		\displaystyle
		s^2 (\alpha - 4) + \left( s^2\alpha - \frac{\beta^2 (\alpha-\eta^2)}{1 - \omega^2} \right) T_s(\omega) + \left( 2\alpha + \frac{\omega \beta (\alpha-\eta^2)}{s^2 (1 - \omega^2)} \right) T_s'(\omega) = 0.
	\end{cases}
\end{equation}
This system can be solved by any root-finding algorithm (for example, Newton's method). Vector $\left(\eta, 1 + \varepsilon / s^2, 1 + \varepsilon \right)^T$ can be used as initial value. For example, solution of this system (with double-precision floating-point numbers) for $s = 5, \eta = 0.95$ is the vector
\begin{equation}\label{eq:examplesol}
	\begin{pmatrix}
		0.950022296412323 \\
		1.0020498847775692 \\
		1.053083013172171
	\end{pmatrix}.
\end{equation}
The resulting polynomials have the form
\begin{equation}\label{eq:exampleR1R0}
	\begin{aligned}
		R_s^1(\mu) &= 1.949130847897793 + 1.0169295750648126\mu \\
		&+ 0.17002420291058604\mu^2 + 0.009987615599077876\mu^3 \\
		&+ 0.00023977479170518486\mu^4 + 0.000002015889739363028\mu^5 \\
		R_s^0(\mu) &= -0.949130847897793 - 0.9660604229626043 \mu\\ 
		&- 0.16151920192429445\mu^2 - 0.009488012136354805\mu^3 \\
		&- 0.00022778070612777503\mu^4 - 0.00000191505030634093\mu^5
	\end{aligned}
\end{equation}
(see Figure~\ref{fig2}).

\section{Error constants and stability}
\label{errstab}

Error constants for methods, based on polynomials \eqref{eq:initialR1R0}, can be easily obtained:
\begin{equation}\label{eq:errconst}
	C_s = \frac{8}{6} - \left( \frac{a}{6} + \frac{r_1^1}{2} + r_2^1 + r_3^1 + r_3^0 \right) = \frac{1}{3} + \frac{1}{6s^2}.
\end{equation}

Error constants for methods, based on damped polynomials \eqref{eq:dampedR1R0}, were obtained numerically. They are presented in Table~\ref{tab:errstab}.

For the stability of the method, the following inequalities must be satisfied:
\begin{equation}
	\begin{aligned}
		-\eta^2 T_s\left( \omega + \beta \frac{\mu}{s^2} \right) &< 1, \\
		\alpha \left(1 + T_s\left( \omega + \beta \frac{\mu}{s^2} 	\right) \right) &> -1 - \eta^2 T_s\left( \omega + \beta \frac{\mu}{s^2} \right)
	\end{aligned}
\end{equation}
or
\begin{equation}\label{eq:stabineq}
	T_s\left( \omega + \beta \frac{\mu}{s^2} \right) > -\frac{1}{\eta^2}, \quad T_s\left( \omega + \beta \frac{\mu}{s^2} \right) > -\frac{1 + \alpha}{\alpha + \eta^2}.
\end{equation}
The second inequality in \eqref{eq:stabineq} is stricter for all $\alpha > 0$. So, we need to find a $\mu = -l_s$ such that the second inequality in \eqref{eq:stabineq} becomes equality:
\begin{equation}\label{eq:stabeq}
	-\cosh \left( s \arccosh \left( -\omega + \beta \frac{l_s}{s^2} \right) \right) = -\frac{1 + \alpha}{\alpha + \eta^2}.
\end{equation}
Solving this equation, we get
\begin{equation}\label{eq:stabintlength}
	l_s = \frac{\omega + \cosh \left( \displaystyle \frac{1}{s} \arccosh \left( \frac{1 + \alpha}{\alpha + \eta^2} \right) \right)}{\beta} s^2.
\end{equation}
The lengths of stability intervals for different parameters $s$ are presented in the Table~\ref{tab:errstab}.

\begin{table}
	\centering
	\caption{The error constants and stability parameters for polynomials \eqref{eq:dampedR1R0}, $\varepsilon = 0.05$}
	\begin{tabular}[c]{|llrr|}
		\hline
		Degree & Error constant & Stability interval & Value \\
		$s$ & $C_s$ & length $l_s$ & $l_s / s^2$ \\
		\hline
		2 & 0.36594 & 7.6531 & 1.913275 \\
		5 & 0.32949 & 47.5779 & 1.903115 \\
		10 & 0.324278 & 190.1654 & 1.901654 \\
		20 & 0.322975 & 760.5155 & 1.901289 \\
		50 & 0.32261 & 4752.9663 & 1.901187 \\
		100 & 0.322558 & 19011.7189 & 1.901172 \\
		200 & 0.322545 & 76046.7294 & 1.901168 \\
		500 & 0.322542 & 475291.8031 & 1.901167 \\
		1000 & 0.322541 & 1901167.0661 & 1.901167 \\
		\hline
	\end{tabular}
	\label{tab:errstab}
\end{table}

Comparing the stability regions with \cite{abdulle2}, we find that the stability interval of our methods is about 2.35 times larger than in ROCK2.

\section{Usage of orthogonal polynomials}
\label{ortpolys}

Following the idea of \cite{vdhouwen}, polynomials \eqref{eq:dampedR1R0} can be written as
\begin{equation}\label{eq:orthoR1R0}
	R_s^1(\mu) = \alpha + T_s(\omega) P_s^1(\mu), \quad
	R_s^0(\mu) = T_s(\omega) P_s^0(\mu),
\end{equation}
where
\begin{equation}
	P_s^1(\mu) = \frac{\alpha}{T_s(\omega)} T_s\left( \omega + \frac{\beta}{s^2} \mu \right), \quad
	P_s^0(\mu) = -\frac{\eta^2}{T_s(\omega)} T_s\left( \omega + \frac{\beta}{s^2} \mu \right)
\end{equation}
or
\begin{equation}\label{eq:P1P0}
	\begin{aligned}
		P_0^1(\mu) &= \alpha, \qquad P_1^1(\mu) = \alpha + \frac{\alpha \beta}{\omega s^2} \mu, \\
		P_0^0(\mu) &= -\eta^2, \quad P_1^0(\mu) = -\eta^2 - \frac{\eta^2 \beta}{\omega s^2} \mu, \\
		P_j^i(\mu) &= 2 \left( \omega + \frac{\beta}{s^2} \mu \right) \frac{T_{j-1}(\omega)}{T_j(\omega)} P_{j-1}^i(\mu) - \frac{T_{j-2}(\omega)}{T_j(\omega)} P_{j-2}^i(\mu).
	\end{aligned}
\end{equation}

Consider two-step Runge-Kutta methods of the form
\begin{equation}\label{eq:recmethod}
	\begin{aligned}
		v_0 &= \tilde{a} y_n + (1-\tilde{a}) y_{n-1}, \\
		v_1 &= v_0 + h \tilde{m}_1 f(x_n + c_0 h, v_0), \\
		v_2 &= m_2 v_1 + (1-m_2)v_0 + h \tilde{m}_2 f(x_n + c_1 h, v_1), \\
		v_3 &= m_3 v_2 + (1-m_3)v_1 + h \tilde{m}_3 f(x_n + c_2 h, v_2), \\
		\vdots \\
		v_s &= m_s v_{s-1} + (1-m_s)v_{s-2} + h \tilde{m}_s f(x_n + c_{s-1} h, v_{s-1}), \\
		y_{n+1} &= a y_n + b v_s,
	\end{aligned}
\end{equation}
where
\begin{equation}\label{eq:recc}
	c_0 = \tilde{a} - 1, \quad
	c_1 = \tilde{a} - 1 + \tilde{m}_1, \quad
	c_j = m_j c_{j-1} + (1 - m_j)c_{j-2} + \tilde{m}_j, \,\, j \ge 2.
\end{equation}

Polynomials $R^1$ and $R^0$ of these methods have the form
\begin{equation}\label{eq:recR1R0}
	\begin{aligned}
		R_s^1(\mu) &= a + b \left( \tilde{R}_s^1(\mu) \right), \quad
		R_s^0(\mu) = b \left( \tilde{R}_s^0(\mu) \right), \\
		\tilde{R}_0^1(\mu) &= \tilde{a}, \quad
		\tilde{R}_1^1(\mu) = \tilde{a} + \tilde{m}_1 \tilde{a} \mu, \\
		\tilde{R}_0^0(\mu) &= 1 - \tilde{a}, \quad
		\tilde{R}_1^0(\mu) = (1 - \tilde{a}) + \tilde{m}_1 (1 - \tilde{a}) \mu, \\
		\tilde{R}_j^i(\mu) &= \left( m_j + \tilde{m}_j \mu \right) \tilde{R}_{j-1}^i(\mu) + \left( 1 - m_j \right) \tilde{R}_{j-2}^i(\mu)
	\end{aligned}
\end{equation}
(see \cite{vdhouwen} for comparison).

Comparing \eqref{eq:recR1R0} with \eqref{eq:orthoR1R0}-\eqref{eq:P1P0} we obtain
\begin{equation}\label{eq:recparams}
	\begin{aligned}
		a = \alpha, \quad 
		b = \left( \alpha - \eta^2 \right) T_s(\omega)&, \quad
		\tilde{a} = \frac{\alpha}{\alpha - \eta^2}, \quad
		\tilde{m}_1 = \frac{\beta}{\omega s^2}, \\
		m_j = 2 \omega \frac{T_{j-1}(\omega)}{T_j(\omega)}, \quad
		\tilde{m}_j &= 2 \frac{\beta}{s^2} \frac{T_{j-1}(\omega)}{T_j(\omega)}, \quad
		j \ge 2.
	\end{aligned}
\end{equation}

For example, for polynomials \eqref{eq:exampleR1R0} we construct the method \eqref{eq:recmethod} with
\begin{equation}
	\begin{tabular}{ll}\label{eq:exampleparams}
		$\begin{aligned}
			\tilde{a} &= 19.991085619464535, \\
			a &= 0.950022296412323, \\
			b &= 0.04997770358767691,
		\end{aligned}$
		& 
		$m = 
		\begin{pmatrix}
			1.9918588786954916 \\ 
			1.9838492426656018 \\
			1.9760315849167438 \\
			1.9684604922450784
		\end{pmatrix},$
		\\ \\
		$\tilde{m} = 
		\begin{pmatrix}
			0.04203714921461939 \\ 
			0.08373206889818684 \\
			0.08339536663324355 \\
			0.08306673458794599 \\
			0.08274846743558949
		\end{pmatrix},$
		& 
		$c = 
		\begin{pmatrix}
			18.991085619464535 \\ 
			19.033122768679153 \\
			19.158549757260907 \\
			19.365346371620134 \\
			19.65025313347653
		\end{pmatrix}.$
		\\
	\end{tabular}
\end{equation}

\ref{code} contains Wolfram Mathematica’s code for obtaining methods \eqref{eq:recmethod} for given parameters $s, \varepsilon$.

\section{Numerical experiments}
\label{numexps}

In all our experiments we use constant step size and reference solutions computed by Wolfram Mathematica’s \texttt{NDSolve}. An additional starting point $y_1$ were taken from this reference solution. For each method we perform a series of constant-step integrations with decreasing step size $h$ and calculate the maximum norm of the error at the endpoint.

Our constant-step integrations can't handle sudden changes in solution components correctly. Therefore, we will choose intervals of integration that do not contain such singularities.

We chose the following stiff problems:
\begin{enumerate}
	\item VDPOL \cite[p. 144]{hairer1}. This problem contains sharp change in the second component of the solution near points $0, 0.8, 1.6, \dots$. So, we start our two-step method from point $x = 0.1$ and finish it at point $x_{out} = 0.6$.
	\item ROBER \cite[p. 144]{hairer1}. Solution of this problem changing more and more slowly and we can take a large integration segment. We start from point $x = 1000$ and finish at point $x_{out} = 2000$.
	\item HIRES \cite[pp. 144-145]{hairer1}. The components of the solution of this problem do not have sharp changes on the segment $[20, 270]$. We will take it to test our method.
	\item Burgers' equation \cite{abdulle4}. We took parameter $\mu = 0.005$ as in \cite{repnikov} but leave the integration interval unchanged: $[0, 2.5]$.
\end{enumerate}

\begin{figure}
	\includegraphics[width=\textwidth]{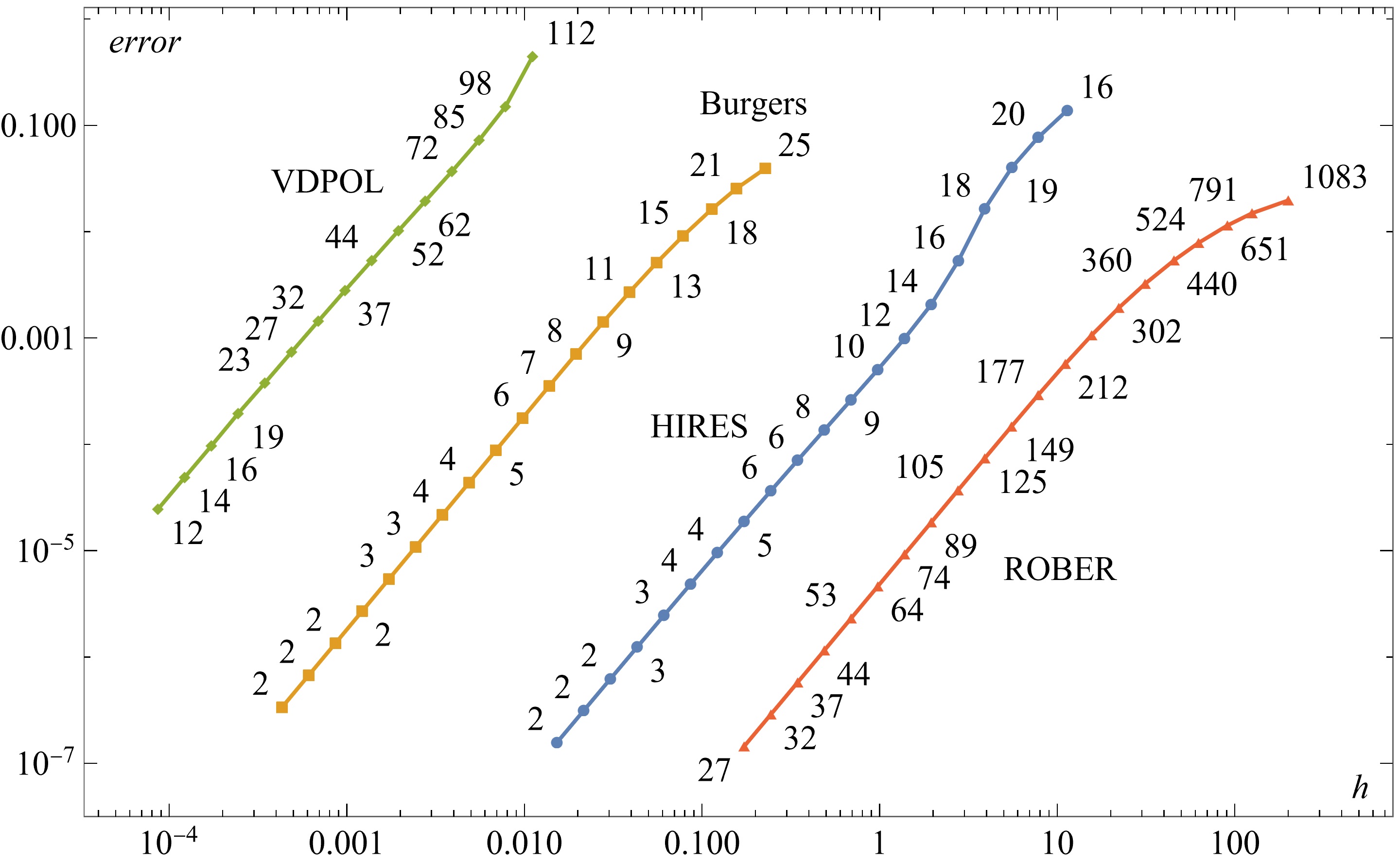}
	\caption{Accuracy and stability of methods \eqref{eq:recmethod}}
	\label{fig3}
\end{figure}

As we can see from \eqref{eq:errconst} and Table~\ref{errstab}, error constants are almost independent of the number of stages (as in the case of one-step methods). This is also confirmed in practice: in the case of stability, the results obtained with different parameters are almost identical. Therefore, we are interested in the number of stages required to achieve stability at the current step size $h$, see Figure~\ref{fig3}. Point labels are equal to the minimum number of stages required to achieve stability at a given step size. For example, to integrate Burgers’ equation with $h = 0.078125$, at least 15 stages of method \eqref{eq:recmethod} are needed. Every method with $\ge 15$ stages gives almost the same result with error $\approx 0.0091$.

It can be seen from the chart that each halving of the step size entails a refinement of the solution by about 4 times which is fully consistent with the second order method. It is also seen that doubling the number of internal stages allows increasing the integration step by 4 times.

\section{Conclusion}
\label{concl}

We have presented new second order two-step Runge-Kutta methods with extended stability interval. In fact, theory of this methods is very close to widely known one-step stabilized methods. The numerical experiments asserted the theoretical properties of accuracy and stability of the constructed methods.

Of course, they are not yet suitable for practical calculations. However, their potential prospects are obvious.

\appendix

\section{Mathematica code for computing method \eqref{eq:recmethod} parameters}
\label{code}

\begin{lstlisting}
getMethod[s_, $\varepsilon$_] := Block[{$\eta$, c, $\alpha$, $\omega$, $\beta$},
   $\eta$ = 1 - $\varepsilon$;
   $R_1$ = $\alpha$ (1 + Cosh[s ArcCosh[$\omega$ + $\beta$ (x / s^2)]]);
   $R_0$ = -$\eta$^2 Cosh[s ArcCosh[$\omega$ + $\beta$ (x / s^2)]];
   $r_1$ = Table[D[$R_1$, {x,j}] / j! /. Rule[x,0], {j,0,3}];
   $r_0$ = Table[D[$R_0$, {x,j}] / j! /. Rule[x,0], {j,0,3}];
   oC = {
      $r_1$[[1]] + $r_0$[[1]] == 1,
      $r_1$[[1]] + $r_1$[[2]] + $r_0$[[2]] == 2,
      $r_1$[[1]]/2 + $r_1$[[2]] + $r_1$[[3]] + $r_0$[[3]] == 2
   };
   vars = N[FindRoot[
      oC, {{$\alpha$, $\eta$}, {$\omega$, 1 + $\varepsilon$ / s^2}, {$\beta$, 1 + $\varepsilon$}}, 
      WorkingPrecision -> 2 MachinePrecision
   ]];
   Print[vars];
   $\alpha$ = $\alpha$ /. vars;
   $\omega$ = $\omega$ /. vars;
   $\beta$ = $\beta$ /. vars;
   len = s^2 (Cosh[1/s ArcCosh[(1 + $\alpha$)/($\alpha$ + $\eta$^2)]] + $\omega$) / $\beta$;
   err = 8/6 - ($r_1$[[1]]/6 + $r_1$[[2]]/2 + $r_1$[[3]] + $r_1$[[4]] + $r_0$[[4]]);
   a = $\alpha$;
   ta = $\alpha$ / ($\alpha$ - $\eta$^2);
   b = ($\alpha$ - $\eta$^2) ChebyshevT[s, $\omega$];
   tm = Join[{$\beta$ / ($\omega$ s^2)}, Table[
      2($\beta$ ChebyshevT[j-1, $\omega$]) / (s^2 ChebyshevT[j, $\omega$]),
      {j,2,s}]];
   m = Table[
      2 $\omega$ ChebyshevT[j-1, $\omega$] / ChebyshevT[j, $\omega$],
      {j,2,s}];
   c = Quiet[RecurrenceTable[{
      c[0] == N[ta-1], c[1] == N[ta-1+tm[[1]]], 
      c[j] == m[[j-1]] c[j-1] + (1 - m[[j-1]]) c[j-2] + tm[[j]]
   }, c, {j, 0, s-1}]];
   Association[Rule["$s$",s],
      Rule["$k$",2], Rule["$order$",2],
      Rule["$a$",a], Rule["$\tilde{a}$",ta], Rule["$b$",b], 
      Rule["$\tilde{m}$",tm], Rule["$m$",m], Rule["$c$",c],
      Rule["$len$",len], Rule["$err$",err],
      Rule["$type$","MultistepRK"]]
];
\end{lstlisting}

\bibliographystyle{elsarticle-num} 
\bibliography{moisa-references}

\end{document}